\documentclass[11pt]{amsart}

\usepackage{amssymb,amscd,amsthm,amsxtra}
\usepackage{latexsym}
\usepackage{epsfig}

\def\be#1{\begin{equation}\label{#1}}

\newtheorem{thm}{Theorem}[section]
\newtheorem{cor}[thm]{Corollary}
\newtheorem{lem}[thm]{Lemma}

\newtheorem{definition}[thm]{Definition}
\theoremstyle{definition}
\newcommand{\comment}[1]{}

\newcommand{\ds}{\mathrm{d}_{\mathrm{s}}}
\newcommand{\dw}{\mathrm{d}_{\mathrm{w}}}
\newcommand{\Rc}{\widetilde{R}}
\newcommand{\Dc}{\mathcal{E}}
\newcommand{\As}{\mathcal{A}_{\mathrm{s}}}
\newcommand{\Aw}{\mathcal{A}_{\mathrm{w}}}
\newcommand{\Xw}{X_{\mathrm{w}}}
\newcommand{\Hw}{H_{\mathrm{w}}}
\newcommand{\ddt}{\frac{d}{dt}}
\newtheorem{defn}[thm]{Definition}
\theoremstyle{remark}

\numberwithin{equation}{section}

\begin{document}

\title[Global attractor of a dyadic model] {An inviscid
dyadic model of turbulence: the global attractor}

\author{Alexey Cheskidov}
\address{Department of Mathematics,
University of Michigan, Ann Arbor, MI  48109-1043}
\email{\tt acheskid@umich.edu}

\author{Susan Friedlander}
\address{Department of Mathematics, Statistics, and
Computer Science, Chicago, IL 60607-7045}
\email{\tt susan@math.uic.edu}

\author{Nata\v{s}a Pavlovi\'{c}}
\address{Department of Mathematics, 
Princeton University, Princeton, NJ 08544-1000}
\email{\tt natasa@math.princeton.edu}

\date{\today}

\begin{abstract}
Properties of an infinite system of nonlinearly coupled
ordinary differential equations are discussed. This
system models some properties present in the equations 
of motion for an inviscid fluid such as the skew symmetry
and the 3-dimensional scaling of the quadratic nonlinearity. 
In a companion paper \cite{CFP} it is proved that 
every solution for the system with forcing
blows up in finite time in the Sobolev $H^{5/6}$ norm. 
In this present paper, it is proved that after the
blow-up time all solutions stay in $H^s$, $s < 5/6$ 
for almost all time and the energy dissipates. 

Moreover, it is proved that the unique equilibrium is an exponential global attractor.
\end{abstract}

\maketitle

\setcounter{equation}{0}
\section{Introduction} 

One of the outstanding open questions in fluid dynamics is existence and
uniqueness of solutions to the Cauchy problem for the three-dimensional 
Euler equations
\begin{equation} \label{Euler}
\frac{\partial u}{\partial t} =-(u \cdot \nabla)u - \nabla p, \qquad \nabla \cdot u=0.
\end{equation}

In last few decades simplified models that capture
some properties of fluid equations have been proposed and studied. In this
article we analyze one of these models, a so called ``dyadic'' model for the
equations of  fluid motion. We study the following dyadic model: 
\begin{equation} \label{Bsystem}
\begin{split}
\frac{d a_j}{d t} & = 2^{\frac{5(j-1)}{2}} a_{j-1}^{2} 
- 2^{\frac{5j}{2}} a_{j} a_{j+1} +f_j, \; \; j >0 \\
\frac{d a_0}{d t} & = - a_0 a_1 + f_0, 
\end{split}
\end{equation}
where the force $f$ is chosen so that 
$f_0 > 0$ and $f_j=0$ for $j>0$ for simplicity.

The model \eqref{Bsystem} without forcing is a special case of the infinite dimensional 
dynamical system 
\be{1dyad} 
\frac{d a_j}{d t} = \lambda^{j-1} a_{j-1}^{2} - \lambda^j a_{j} a_{j+1}.
\end{equation} 
Such an inviscid model has been studied recently in a number of articles
including \cite{FP,KP,KZ,W}. Variants of the model that include viscosity are discussed in
\cite{C,KPalpha}. Analysis of more general ``shell" models
and motivation in terms of turbulence modeling can be found in  
\cite{BJPV, CLT, DS,Gledzer,O,OY}.

In a companion paper \cite{CFP} we presented a motivation for the model 
\eqref{Bsystem} from the Fourier space Euler equations (1.1) in 3-dimensions. 
The coefficient $a_j^2(t)$ is the total energy in the frequency space shell 
$2^j \leq |k| < 2^{j+1}$. In this context $l^2$ and 
$H^s$, respectively the energy and Sobolev norms, are defined as 
\be{norms} 
\|a(t)\|_{l^2} = \left( \sum_{j=0}^{\infty} a_j^2(t) \right)^{1/2}, \; \; \; 
\|a(t)\|_{H^s} = \left( \sum_{j=0}^{\infty} 2^{2sj} a_j^2(t) \right)^{1/2}. 
\end{equation}  
The nonlinear terms on the right hand side of \eqref{Bsystem} 
retain important features of the advective term in the Euler 
equation, namely bilinearity and skew-symmetry. The presence 
of the specific quadratic term $2^{\frac{5(j-1)}{2}} a_{j-1}^2$ ensures 
a certain monotonicity (see also, \cite{BB} and \cite{P})
in the cascade of energy through the scales $j$. 
We defined a regular solution for the model \eqref{Bsystem} to be a 
solution with bounded $H^{5/6}$ norm and such solutions satisfy 
the energy equality. In \cite{CFP} we proved that: 
\begin{enumerate} 
\item[(a)] There exists a unique fixed point to \eqref{Bsystem} 
and the fixed point is not in $H^{5/6}$. 

\item[(b)] Every regular solution approaches the fixed point 
in the $l^2$ norm. 

\item[(c)] Every solution blows up in finite time in the $H^{5/6}$ norm. 
\end{enumerate} 

As Mattingly et al \cite{MS} observe in their recent analysis of an infinite
linear dynamical system, the most interesting features of such models 
belong to solutions {\em{after}} the time of blow-up when some norm becomes infinite. 
This is particularly true in the context of models that illustrate behavior 
that has been proposed to characterize hydrodynamic turbulence in the works of,
for example, Kolmogorov \cite{K}, Onsager \cite{Ons}, Frisch \cite{Fr}, Robert \cite{R}, 
Constantin et al \cite{CET} and  Eyink-Sreenivasan \cite{ES}. 

In our present paper we study the solutions of \eqref{Bsystem} 
{\em{after}} the time of blow-up in $H^{5/6}$. We prove: 
\begin{enumerate} 
\item The $H^{s}$ norms for $s < 5/6$ are locally square integrable in time. 

\item The solutions dissipate energy. 

\item The unique fixed point $\{a_j\} = \{2^{-5j/6} 2^{5/12} \sqrt{f_0} \}$ 
is an exponential global attractor. 

\end{enumerate} 
The existence of a global attractor for an inviscid system is, perhaps, surprising. 
However it is exactly consistent with the concept of anomalous or turbulent dissipation
conjectured by Onsager \cite{Ons}. As we discuss in Section \ref{secOns}, after the 
time of blow up in $H^{5/6}$ the energy spectrum exactly reproduces Kolmogorov's law: 
\be{KL} 
E(|k|) = c_0 \bar{\epsilon}^{2/3} |k|^{-5/3},
\end{equation}  
where $\bar{\epsilon}$ is the average of the energy dissipation rate. 

\subsection*{Notation} 
For notational convenience we adopt $\lambda^j$ as the scaling parameter in 
the analysis performed in sections 2 - 5. We do this to illustrate that 
the results are qualitatively independent of the exact choice
of $\lambda$ (which depends on the spatial dimension and the 
construction of the model). The proofs of results in section \ref{secfix}
require $\lambda < 2^3$. As we discussed in \cite{CFP} 
the relevant $\lambda$ for the 3-dimensional model is 
$2^{5/2}$. This exponent determines the values of the
exponent of the fixed point and the critical Sobolev space 
exponent $H^{5/6}$. The exponent also reproduces the exponents in the
Kolmogorov's law \eqref{KL}.

\subsection*{Acknowledgements} The authors would like to thank 
Marie Farge, Jonathan Mattingly, Kai Schneider and Eric Vanden-Eijnden
for very helpful discussions. S.F. was partially supported by NSF grant 
number DMS 0503768. N.P. was partially supported by NSF grant 
number DMS 0304594.

%%%%%%%%%%%%%%%%%%%%%%%%%%%%%%%%
\section{Functional setting} \label{s:setting}
%%%%%%%%%%%%%%%%%%%%%%%%%%%%%%%%
Let us denote $H=l^2$ with the usual scalar product and norm:
\begin{equation}
(a,b):= \sum_{j=0}^{\infty} a_jb_j, \qquad |a|:=\sqrt{(a,b)}.
\end{equation}
The norm $|a|$ will be called the energy norm.
%Let $A:D(A) \to H$ be the Laplace operator defined by
%\begin{equation}
%(Aa)_n = \lambda^{2\alpha n} u_n, \qquad n\geq 1,
%\end{equation}
%for some $\lambda >1$.
%The domain $D(A)$ of this operator is a dense subset of $H$.
%Note that $A$ is a positive definite operator
%whose eigenvalues are
%\begin{equation}
%0<\lambda^{2\alpha} \leq \lambda^{4\alpha} \leq \lambda^{6\alpha} \leq \dots
%\end{equation}
%Let $H^{s}=A^{-\gamma/(2\alpha)}H$ endowed with the
%following scalar product and norm: 
Let 
\begin{equation}
((a,b))_{s}:= \sum_{j=0}^{\infty} 2^{2s j}a_jb_j, \qquad \|u\|_{s}:=\sqrt{((u,u))_s}.
\end{equation}
%In the special case $\gamma=\alpha$, let $V=H^{\alpha} =A^{-1/2}H$ and
%\begin{equation}
%((u,v)):=((u,v))_{\alpha}, \qquad \|u\|:=\|u\|_\alpha.
%\end{equation}
%This double norm $\|u\|$ will be called the enstrophy norm.
%Note that we have an equivalent of the Poincar\'{e} inequality
%\begin{equation}
%|u|^2 \leq \frac{1}{\lambda^{2\alpha}} \|u\|^2.
%\end{equation}
We fix $\lambda=2^{5/2}$ and let
\begin{equation}
\ds(a,b):=|a-b|, \qquad
\dw(a,b):= \sum_{j=0}^\infty \frac{1}{\lambda^{(j^2)}}
\frac{|a_j-b_j|}{1 + |a_j-b_j|},
\qquad a,b \in H.
\end{equation}
Here, $\ds$ is a strong distance, and $\dw$ is a weak distance
that induces a weak topology on any bounded subset of $H$. Hence,
a bounded sequence $\{a^k\} \subset H$ converges to $a\in H$ weakly,
i.e.,
\begin{equation}
\lim_{k\to \infty} (a^k, b) = (a,b),  \qquad \forall b \in H,
\end{equation}
if and only if
\begin{equation}
\dw(a^k, a) \to 0 \qquad \text{as} \qquad  k\to \infty.
\end{equation}
%We also recall that if $u^k \to u$ weakly in $H$ as $k\to \infty$, then
%\begin{equation}
%\liminf_{k \to \infty}|u^k| \geq |u|.
%\end{equation}
Let
\begin{equation}
C([0,T];\Hw):=\{a(\cdot): [0,T] \to H, a_j(t) \mbox{ is continuous for all }
j\}
\end{equation}
endowed with the distance
\begin{equation}
d_{C([0, T];\Hw)}(a,b) = \sup_{t\in[0,T]}\dw(a(t),b(t)). 
\end{equation}
Let also
\begin{equation}
C([0,\infty);\Hw):=\{a(\cdot): [0,\infty) \to H, a_j(t) \mbox{ is continuous for all }
j\}
\end{equation}
endowed with the distance
\begin{equation}
d_{C([0, \infty);\Hw)}(a,b) = \sum_{T\in \mathbb{N}} \frac{1}{\lambda^T} \frac{\sup\{\dw(a(t),b(t)):0\leq t\leq T\}}
{1+\sup\{\dw(a(t),b(t)):0\leq t\leq T\}}.
\end{equation}

%%%%%%%%%%%%%%%%%%%%%%%%%%%%%%%%%%%%%%%%%%%%%%%%%%%%%%%%%%%%%%

\section{Weak solutions}

\begin{definition}
A weak solution on $[0,T]$ (or $(0, \infty)$, if
$T=\infty$) of \eqref{Bsystem} is an $H$-valued
function $a(t)$ defined for $t \in [0, T]$, such that $a_j \in C^1([0,T])$
and $a_j(t)$ satisfies \eqref{Bsystem}  for all $j$.
\end{definition}

Note that since the nonlinear term has a finite number of terms, the notions of
a weak solution and a classical solution (of a system of ODEs) coincide.
Hence, the weak solutions will be called solutions in the remainder of the paper. 
Note that if $a(t)$ is a solution on $[T,\infty)$, then automatically
$a_j \in C^{\infty}([T,\infty))$.

\begin{thm}[Global existence] \label{thm:Leray}
For every $a^0 \in H$,
there exists a solution of \eqref{Bsystem} on $[0,\infty)$ with $a(0)=a^0$.
\end{thm}
\begin{proof}
Let $u^0 \in H$ and $T>0$ be arbitrary.
We will show the existence of a solution on $[0,T]$ by taking a limit of the Galerkin
approximation $a^k(t)=(a^k_0(t),\dots,a^k_k(t),0,0,\dots)$ with
$a^k_j(0)=a^0_j$ for $j=1,2,...,k$, which satisfies
\begin{equation} \label{galerkin}
\left\{
\begin{aligned}
&\ddt a^k_j  - \lambda^{j-1} (a^k_{j-1})^2 +
\lambda^{j} a^k_j a^k_{j+1}
=f_j, \qquad j\leq k-1,\\
&\ddt a^k_k  - \lambda^{k-1} (a^k_{k-1})^2 =f_k,
%&u_n(0) = u_n^0.
\end{aligned}
\right.
\end{equation}
where $a^j_{-1}=0$ and $\lambda=2^{5/2}$.
From the theory of ordinary differential equations we know that 
there exists a unique solution $a^k(t)$ to \eqref{galerkin} on $[0,T]$.
%Consider a sequence of Galerkin approximations $\{u^N\}$.
We will show that a sequence of the Galerkin approximations $\{a^k\}$ is
weakly equicontinuous. Indeed, it is clear that there exists $M$, such that
\begin{equation}
a^k_j(t) \leq M, \qquad \mbox{ for all } t\in[0,T] \mbox{ and all } j, k. 
\end{equation}
Therefore,
\begin{equation}
\begin{split}
|a^k_j(t) - a^k_j(s)| &\leq \left| \int_s^t \left(
\lambda^{j-1}(a^k_{j-1}(\tau))^2 -\lambda^j a^k_j(\tau) a^k_{j+1}(\tau) + f_j\right) \, d\tau\right|\\
&\leq (\lambda^{j-1} M^2 + \lambda^jM^2 + f_j)
|t-s|,
\end{split}
\end{equation}
for all $j$, $k$ and all $0\leq t \leq s \leq T$. Thus,
\begin{equation}
\begin{split}
\dw(a^k(t), a^k(s))&= \sum_{j=0}^\infty \frac{1}{\lambda^{(j^2)}}
\frac{|a^k_j(t)-a^k_j(s)|}{1 + |a^k_j(t)-a^k_j(s)|}\\
%&\leq \sum_{n=1}^N \lambda^{-2n} \left| \int_t^s 2M^2 \lambda^{n+1}  d\tau \right|
&\leq c|t-s|,
\end{split}
\end{equation}
for some constant $c$ independent of $k$. Hence, $\{a^k\}$ is an equicontinuous sequence
of functions in $C([0,T];\Hw)$ with bounded initial data. Therefore,
the Ascoli-Arzela theorem implies that
$\{a^k\}$ is relatively compact in $C([0,T]; \Hw)$.
Hence, passing to a subsequence,
we obtain that there exists a weakly continuous $H$-valued
function $a(t)$, such that
\begin{equation} \label{e:weakconv}
a^{k_n} \to a \qquad \mbox{as} \qquad k_n \to \infty  \qquad \mbox{in}
\qquad  C([0,T]; \Hw).
\end{equation}
In particular, $a^{k_n}_j(t) \to a_j(t)$ as $k_n \to \infty$, for all $j, t \in [0, T]$.
Thus, $a(0) = a^0$.
In addition, note that
\begin{equation}
a^{k_n}_j(t) =a^{k_n}_j(0) + \int_0^t (
\lambda^{j-1} (a^{k_n}_{j-1})^2 - \lambda^j a^{k_n}_j a^{k_n}_{j+1} +f_j) \, d\tau,
\end{equation}
for $j \leq k_n-1$. Taking the limit as $k_n \to \infty$, we obtain
\begin{equation}
a_j(t) = a_j(0) + \int_0^t (
\lambda^{j-1} a_{j-1}^2 - \lambda^j a_j a_{j+1} +f_j) \, d\tau.
\end{equation}
Since $a_j(t)$ is continuous, it follows that
$a_j \in C^1([0,T])$ and satisfies \eqref{Bsystem}.
\end{proof}

\begin{thm}[Energy inequality] \label{t:iquality}
Let $a(t)$ be a solution of \eqref{Bsystem} with $a_j(0)\geq0$.
Then $a_j(t) > 0$ for all $t>0$, and 
$a(t)$ satisfies the energy inequality
\begin{equation} \label{ee}
|a(t)|^2 \leq |a(t_0)|^2 + 2\int_{t_0}^t (f, a(\tau)) \, d\tau,
\end{equation}
for all $0 \leq t_0 \leq t$.
\end{thm}
\begin{proof}
A general solution of \eqref{Bsystem} can be written as
\begin{multline} \label{e:gensol}
a_j(t)=a_j(0)\exp\left(-
\int_{0}^t \lambda^{j} a_{j+1}(\tau) \, d\tau\right)\\
+ \int_{0}^t\exp\left(-\int_{s}^{t} \lambda^{j}
a_{j+1}(\tau) \, d\tau\right) (f_j+\lambda^{j-1} a_{j-1}^2(s) )\, ds,
\end{multline}
where $\lambda=2^{5/2}$.
Recall that $f_j \geq 0$ for all $j$.
Since $a_j(0) \geq 0$ for all $j$, then  $a_j(t)\geq 0$ for all $j$, $t>0$.
Moreover, since $f_0>0$, we have $a_0(t)>0$ for all $t>0$ and, consequently,
$a_j(t)>0$ for all $j$, $t>0$.
Hence, multiplying \eqref{Bsystem} by $a_j$,
taking a sum from $0$ to $N$, and integrating between $t_0$ and $t$,
we obtain
\begin{equation} \label{e:eipar}
\begin{split}
\sum_{j=0}^N a_j(t)^2  - \sum_{j=0}^N a_j(t_0)^2
 &= - 2\int_{t_0}^t  \lambda^{N} a_N^2
a_{N+1} \, d\tau  +2\int_{t_0}^t\sum_{j=0}^N  f_j a_j \, d\tau\\
&\leq 2\int_{t_0}^t\sum_{j=0}^N f_j a_j \, d\tau.
\end{split}
\end{equation}
Taking the limit as $N \to \infty$, we obtain \eqref{ee}.

\end{proof}

%%%%%%%%%%%%%%%%%%%%%%%%%%%%%%%%%%%%%%%%%%%%%%%%%%%%%%%%%%%%%%%%%%%%%%%%%%

\section{Fixed point} \label{secfix}

%Here we present a proof of Theorem \ref{nstab}. 
%We consider the fully nonlinear system %\eqref{bj} 
Given a solution $a(t)$ of \eqref{Bsystem} with arbitrary initial data $a_j(0) \geq0$,
let
\begin{equation}
b_j(t):=a_j(t) - \lambda^{-j/3} , \qquad j\geq 0,
\end{equation}
where $\lambda=2^{5/2}$. Let also
\begin{equation} \label{defb}
\begin{split} 
d_j&=\lambda^{\frac{j-1}{3}}(\lambda^{\frac{1}{6}} b_{j} - \lambda^{-\frac{1}{6}} b_{j-1}), \qquad j\geq 1,\\
d_0&= \lambda^{-\frac{1}{6}} b_0.
\end{split}
\end{equation}
%Now define the following function
%\begin{equation} \label{defphi}
%\phi(t)=
%\begin{cases}
%\frac{1}{4}\|d(t)\|_{l^2}^2, & \text{if } \|d(t)\|_{l^2}^2 < \infty,\\
%1, & \text{if } \|d(t)\|_{l^2}^2 = \infty.
%\end{cases}
%\end{equation}
%Note that
%\begin{equation}
%\lambda^{-\frac{j-1}{3}}d_j= \lambda^{\frac{1}{6}} b_j-\lambda^{-\frac{1}{6}}b_{j-1}.
%\end{equation}
%Therefore, we have
%\begin{equation}
%\begin{split}
%\lambda^{\frac{1}{3}}d_0&= \lambda^{\frac{1}{6}} b_0,\\
%d_1&= \lambda^{\frac{1}{6}} b_1 -  \lambda^{-\frac{1}{6}} b_0,\\
%\lambda^{-\frac{2-1}{3}}d_2&= \lambda^{\frac{1}{6}} b_2 - \lambda^{-\frac{1}{6}} b_1,\\
%& \vdots \\
%\lambda^{-\frac{j-1}{3}}d_j&= \lambda^{\frac{1}{6}} b_j -  \lambda^{-\frac{1}{6}} b_{j-1}.
%\end{split}
%\end{equation}
We sum the expressions for $d_l$, with $0 \leq l \leq j$ and thanks to \eqref{defb} obtain: 
\begin{equation} \label{eq:bd}
b_j = \lambda^{\frac{1}{6}-\frac{j}{3}}(d_0+d_1+\dots+d_j).
\end{equation}

\begin{lem}\label{stabilitylemma}
For every $t\geq 0$ we have
\begin{equation}
\frac{d}{dt} \sum_{j=0}^{k} b_{j}(t)^{2} \leq -\sum_{j=0}^{k} d_j(t)^2 + 
d_{k+1}(t)^2,
%\qquad k \geq 1,
\end{equation}
and
\begin{equation}
\frac{d}{dt}\left( \sum_{j=0}^{k-1} b_{j}(t)^{2} +\frac{1}{2}b_k(t)^2 \right)
\leq -\sum_{j=0}^{k-1} d_{j}(t)^{2} + \lambda^{\frac{k-1}{3}} |b(t)|,
\end{equation}
for all $k >1$.
\end{lem}
\begin{proof}
Note that $b(t)$ satisfies the following system of equations:
\begin{equation}
\begin{split} \label{rpbj}
\frac{d b_j}{dt} & = \lambda^{\frac{2j}{3}} 
( 2\lambda^{-\frac{2}{3}} b_{j-1} - \lambda^{-\frac{1}{3}} b_j - b_{j+1})
 +  ( \lambda^{j-1} b_{j-1}^2 - \lambda^j b_j b_{j+1} ), \; \; j\geq 1, \\
\frac{d b_0}{dt} & = -\lambda^{-\frac{1}{3}} b_0 - b_1 -  b_0 b_1.
\end{split} 
\end{equation}
%The energy of the system is given by 
Multiplying it by $b_j$ and taking a sum from $j=0$ to $j=k$ we obtain
\begin{multline} \label{beftel-enbj}
 \frac{1}{2} \frac{d}{dt} \sum_{j=0}^{k} b_j^2 
= -\lambda^{-\frac{1}{3}} b_0^2 - b_0b_1 - b_0^2 b_1 \\
 + \sum_{j=1}^{k} \left( 
\lambda^{\frac{2j}{3}} ( 2\lambda^{-\frac{2}{3}} b_{j-1}b_j - \lambda^{-\frac{1}{3}}b_j^2 - b_jb_{j+1}) 
+   ( \lambda^{j-1} b_{j-1}^2b_j - \lambda^j b_j^2 b_{j+1} ) \right). 
\end{multline}
%We note that for $b_j$ small (i.e. $b_j < \lambda^{-\frac{j}{3}}$) the
%summations in the energy expression \eqref{beftel-enbj} for $b_j$ converge uniformly. 
%Thus the telescoping sum is justified,
%We note that in the nonlinear term $O(\epsilon)$ all the terms
%cancel except the last one, obtaining
It now follows that
\begin{multline} \label{enbj}
\frac{1}{2} \frac{d}{dt} \sum_{j=0}^{k} b_{j}^{2}   
=  -\left(\lambda^{-\frac{1}{3}} b_0^2 + b_0b_1\right) \\
+ \sum_{j=1}^{k} \lambda^{\frac{2j}{3}} 
\left( 2\lambda^{-\frac{2}{3}} b_{j-1}b_j - \lambda^{-\frac{1}{3}}b_j^2 - b_jb_{j+1}\right)
- \lambda^k b_k^2b_{k+1}.
\end{multline} 
Also we can rewrite \eqref{enbj} as 
\be{enbjrew} 
\frac{1}{2} \frac{d}{dt} \sum_{j=0}^{k} b_{j}^{2} 
= - \lambda^{-\frac{1}{3}} \sum_{j=0}^{k} \lambda^{\frac{2j}{3}} b_j^2 + 
\sum_{j=0}^{k} \lambda^{\frac{2j}{3}} b_j b_{j+1}
-2\lambda^{\frac{2k}{3}} b_kb_{k+1}
 - \lambda^k b_k^2b_{k+1}. 
\end{equation} 
However,
\begin{align*} 
& \sum_{j=0}^{k} \lambda^{\frac{2j}{3}} 
(\lambda^{-\frac{1}{6}} b_j - \lambda^{\frac{1}{6}} b_{j+1})^2 \\
&\quad = \sum_{j=0}^{k} \lambda^{\frac{2j}{3}} 
(\lambda^{-\frac{1}{3}} b_j^2 + \lambda^{\frac{1}{3}} b_{j+1}^2 - 2 b_j b_{j+1}) \\
&\quad = -2 \left[ - \lambda^{-\frac{1}{3}} \sum_{j=0}^{k} \lambda^{\frac{2j}{3}} b_j^2   
+ \sum_{j=0}^{k} \lambda^{\frac{2j}{3}} b_j b_{j+1} \right]
- \lambda^{-\frac{1}{3}} b_0^2
+\lambda^{\frac{2k}{3}+\frac{1}{3}}b_{k+1}^2.  
\end{align*} 
Hence \eqref{enbjrew} gives 
\begin{multline} \label{ebjrew1}
\frac{1}{2} \frac{d}{dt} \sum_{j=0}^{k} b_{j}^{2} 
= -\frac{1}{2} \left[
\sum_{j=0}^{k} \lambda^{\frac{2j}{3}} 
(\lambda^{-\frac{1}{6}} b_j - \lambda^{\frac{1}{6}} b_{j+1})^2 
+ \lambda^{-\frac{1}{3}} b_0^2
\right]\\
+\frac{1}{2} \lambda^{\frac{2k}{3} + \frac{1}{3}}b_{k+1}^2
-2\lambda^{\frac{2k}{3}} b_kb_{k+1}
 - \lambda^k b_k^2b_{k+1}.
\end{multline}

%We claim that for every $t\geq 0$ there exists $N$, such that
%\begin{equation}
%\frac{d}{dt} \sum_{j=0}^{k} b_{j}(t)^{2} < 0, \qquad \forall k >N.
%\end{equation}
Now note that since the initial condition of the solution satisfies
$a_j(0) \geq0$ for all $j$, we have that
$a_j(t)> 0$ for all $j$ and $t>0$. Therefore,
\begin{equation} \label{e:temp}
 b_j(t) \geq -\lambda^{-\frac{j}{3}}, \qquad \forall j, t\geq 0.
\end{equation}
Thus,
\begin{equation}
 - \lambda^k b_k^2b_{k+1} \leq \lambda^{\frac{2k}{3}-\frac{1}{3}}b_k^2,
\end{equation}
for all $k$, $t\geq0$.

Hence,
\begin{equation}  \label{e:ddtbj}
\begin{split} 
\frac{1}{2} \frac{d}{dt} \sum_{j=0}^{k} b_{j}^{2} 
& \leq -\frac{1}{2}\sum_{j=0}^{k+1} d_{j}^{2} 
+ \frac{1}{2} \lambda^{\frac{2k}{3} + \frac{1}{3}}b_{k+1}^2 -2\lambda^{\frac{2k}{3}} b_kb_{k+1} +  \lambda^{\frac{2k}{3}-\frac{1}{3}}b_k^2  \\
& = -\frac{1}{2}\sum_{j=0}^{k+1} d_{j}^{2} +  d_{k+1}^2 - \frac{1}{2} \lambda^{\frac{2k}{3}+ \frac{1}{3}}b_{k+1}^2.
%& \leq -\frac{1}{2}\sum_{j=0}^{k+1} d_{j}^{2} +  d_{k+1}^2.
\end{split} 
\end{equation} 
Note that from \eqref{defb} it follows that
\begin{equation} \label{cons-defb}
b_{k+1} = \lambda^{-\frac{k}{3}-\frac{1}{6}}d_{k+1} + \lambda^{-\frac{1}{3}}b_k.
\end{equation}
Now we rewrite \eqref{e:ddtbj} using \eqref{cons-defb} and \eqref{defb} as follows: 
%Thanks to \eqref{defb} we obtain
\begin{equation} \label{eq:sadfsd}
\begin{split}
\frac{d}{dt} \sum_{j=0}^{k} b_{j}^{2} 
 &\leq -\sum_{j=0}^{k} d_{j}^{2} -2\lambda^{\frac{k}{3}-\frac{1}{6}}d_{k+1}b_k
-\lambda^{\frac{2k}{3}-\frac{1}{3}}b_k^2\\
&= -\sum_{j=0}^{k} d_{j}^{2} -2\lambda^{\frac{2k}{3}}b_{k+1}b_k + \lambda^{\frac{2k}{3}-\frac{1}{3}}b_k^2\\
&= -\sum_{j=0}^{k-1} d_{j}^{2} -2\lambda^{\frac{2k}{3}}b_{k+1}b_k
 +\lambda^{\frac{2k}{3}-1}b_{k-1}^2 +2\lambda^{\frac{k}{3}-\frac{1}{2}}d_kb_{k-1}\\
&=  -\sum_{j=0}^{k-1} d_{j}^{2} -2\lambda^{\frac{2k}{3}}b_{k+1}b_k
+2\lambda^{\frac{2k}{3}-\frac{2}{3}}b_kb_{k-1} - \lambda^{\frac{2k}{3}-1}b_{k-1}^2.
\end{split}
\end{equation} 
Note that we also have
\begin{equation} \label{eq:waerwerwea}
\frac{d}{dt} \sum_{j=0}^{k-1} b_{j}^{2} \leq
  -\sum_{j=0}^{k-1} d_{j}^{2}
-2\lambda^{\frac{2k}{3}-\frac{2}{3}}b_kb_{k-1} + \lambda^{\frac{2k}{3}-1}b_{k-1}^2.
\end{equation} 
Adding equations \eqref{eq:sadfsd} and \eqref{eq:waerwerwea} we get
\begin{equation} \label{eq:split}
\begin{split}
\frac{d}{dt}\left( \sum_{j=0}^{k-1} b_{j}^{2} +\frac{1}{2}b_k^2 \right) &\leq
  -\sum_{j=0}^{k-1} d_{j}^{2} -\lambda^{\frac{2k}{3}}b_{k+1}b_k\\
&\leq -\sum_{j=0}^{k-1} d_{j}^{2} + \lambda^{\frac{k}{3}-\frac{1}{3}} |b|.
\end{split}
\end{equation} 

\end{proof}

\begin{thm} \label{decreasing-norm}
For every solution $a(t)$ with the initial data $a(0)\in l^2$, $a_j(0) \geq0$,
and every time interval $[t_1, t_2]$, $0\leq t_1 \leq t_2 $, we
have that
\begin{equation}
|b(t_2)|^2 -|b(t_1)|^2 \leq
- \alpha\int_{t_1}^{t_2} |d(t)|^2 \, dt,
\end{equation}
where $\alpha =2-\lambda^{3/8}$.
%the $l^2$ norm  $\| \epsilon b(t)\|_{l^2}$ is monotonically
%decreasing for all $t\geq0$.
\end{thm}
\begin{proof}
Since $\lambda=2^{5/2}$, we have that $\alpha\in (0,1)$.
First assume that exists $N>0$, such that
\begin{equation} \label{eq:tocontr}
\sum_{j=0}^k b_j(t_2)^2 - \sum_{j=0}^k b_j(t_1)^2 >
-\alpha \int_{t_1}^{t_2} \sum_{j=0}^{k} d_j(t)^2 \, dt,
\end{equation}
for all $k\geq N$.
On the other hand, thanks to Lemma~\ref{stabilitylemma}, we have that
\begin{equation}
\sum_{j=0}^k b_j(t_2)^2 - \sum_{j=0}^k b_j(t_1)^2 \leq
-\int_{t_1}^{t_2} \sum_{j=0}^{k} d_j(t)^2 \, dt
+\int_{t_1}^{t_2} d_{k+1}(t)^2 \, dt.
\end{equation}
Hence,
\begin{equation} \label{eq:intdk}
\int_{t_1}^{t_2} d_{k+1}(t)^2 \, dt > (1-\alpha) \int_{t_1}^{t_2} \sum_{j=0}^{k} d_j(t)^2 \, dt,
\end{equation}
for all $k \geq N$.
Let
\begin{equation}
I_j := \int_{t_1}^{t_2} d_j(t)^2 \, dt.
\end{equation}
>From \eqref{eq:intdk} it follows that
\begin{equation}
I_{k+1} > (1-\alpha)\sum_{j=0}^k I_j, \qquad \forall k\geq N.
\end{equation}
Hence,
\begin{equation}
\begin{split}
I_{N+1} &>0,\\
I_{N+2} &> (1-\alpha) I_{N+1},\\
I_{N+3} &> (1-\alpha)(I_{N+2}+I_{N+1}) > (1-\alpha) (2-\alpha) I_{N+1},\\
&\vdots \\
I_{N+j} & > (1-\alpha)(2-\alpha)^{j-2} I_{N+1}.
\end{split}
\end{equation}
Therefore,
\begin{equation} \label{EQ:jkoljkljkl}
I_{N+j} > \lambda^{\frac{3j}{8}} \frac{1-\alpha}{(2-\alpha)^2}I_{N+1} , \qquad \forall j\geq 2.
\end{equation}
Now let $M=(t_2-t_1)\sup_{t\in[t_1, t_2]} \|b(t)\|_{l_2}$.
Lemma~\ref{stabilitylemma} and \eqref{EQ:jkoljkljkl} imply that 
\begin{equation}
\begin{split}
\sum_{j=0}^{N+j} &b_{j}(t_2)^{2} -\sum_{j=0}^{N+j} b_{j}(t_1)^{2}  +\alpha I_{N+j}\\
&\leq \frac{1}{2}b_{N+j+1}(t_1)^2   -\frac{1}{2}b_{N+j+1}(t_2)^2    - (1-\alpha)I_{N+j} + \lambda^{\frac{N+j}{3}}M\\
&\leq \frac{1}{2}b_{N+j+1}(t_1)^2   -\frac{1}{2}b_{N+j+1}(t_2)^2
  - \lambda^{\frac{3j}{8}} \frac{(1-\alpha)^2}{(2-\alpha)^2}I_{N+1}
+ \lambda^{\frac{N+j}{3}}M,
\end{split}
\end{equation}
for all $j \geq 2$. Note that the right hand side goes to $-\infty$ as $j\to \infty$,
contradicting \eqref{eq:tocontr}.

Therefore, we have shown that for any $N>0$ there exists $k>N$, such that
\begin{equation}
\sum_{j=0}^k b_j(t_2)^2 - \sum_{j=0}^k b_j(t_1)^2 \leq
-\alpha \int_{t_1}^{t_2} \sum_{j=0}^{k} d_j(t)^2 \, dt.
\end{equation}
By the definition of a weak solution,  $b(t) \in l^2$ for all time $t$.
Therefore, taking a limit as $N \to \infty$ and using Levi's convergence theorem,
we obtain that $|d(t)|^2$ is locally integrable and
\begin{equation}
|b(t_2)|^2 -|b(t_1)|^2 \leq
- \alpha\int_{t_1}^{t_2} |d(t)|^2 \, dt,
\end{equation}
which concludes the proof.
\end{proof}

\begin{thm} Let $a(t)$ be a solution of \eqref{Bsystem} with the initial data $a(0)\in l^2$, $a_j(0) \geq0$.
Then $\|a(t)\|^2_s$ is locally integrable on $[0,\infty)$ for all $s < 5/6$.
\end{thm}
\begin{proof}
Thanks to \eqref{eq:bd}, we have that
\begin{align} 
 \sum_{j=0}^\infty 2^{2sj} b_j^2 
&= \sum_{j=0}^\infty 2^{2sj} \lambda^{\frac{1}{3}-\frac{2j}{3}}(d_0+d_1+\dots+d_j)^2 \nonumber \\
&\leq \sum_{j=0}^\infty 2^{2sj} \lambda^{\frac{1}{3}-\frac{2j}{3}} (j+1) (d_0^2+d_1^2+\dots+d_j^2) \label{lamineq}  \\
&\leq 2^{\frac{5}{6}} |d|^2  \sum_{j=0}^\infty 2^{2j(s - \frac{5}{6})} (j+1), \nonumber %\label{eq:bvsd}
\end{align}
where to obtain the last line we used $\lambda = 2^{5/2}$. 
Due to Theorem~\ref{decreasing-norm}, $|d(t)|$ is locally integrable.
Therefore, $\|b(t)\|_s^2$ is locally integrable, provided $s<5/6$. Hence,
$\|a(t)\|_s^2$ is locally integrable for $s<5/6$.
\end{proof}

\begin{thm} \label{t:Expon}
Let $a(t)$ be a solution of \eqref{Bsystem} with $a_j(0)\geq 0$. Then
$a(t)$ exponentially converges in $l^2$ to the fixed point as $t \to \infty$. More precisely,
\begin{equation}
|b(t)|^2 \leq |b(0)|^2 e^{-\beta t},
\end{equation}
for some universal constant $\beta>0$.
\end{thm}
\begin{proof}
Note that inequality \eqref{lamineq} with $s=0$ implies that
\begin{equation}
|b|^2 \leq  |d|^2 \sum_{j=0}^\infty \lambda^{\frac{1}{3}-\frac{2j}{3}} (j+1).
\end{equation}
Therefore, Theorem~\ref{decreasing-norm} yields
\begin{equation}
|b(t_2)|^2 -|b(t_1)|^2 \leq
- \beta \int_{t_1}^{t_2} |b(t)|^2 \, dt,
\end{equation}
where
\begin{equation}
\beta=\frac{ 2-\lambda^{\frac{3}{8}} }{ \sum_{j=0}^\infty \lambda^{\frac{1}{3}-\frac{2j}{3}} (j+1) }.
\end{equation}
Using Granwall's inequality, we conclude that
\begin{equation}
|b(t)|^2 \leq |b(0)|^2 e^{-\beta t}.
\end{equation}
\end{proof}

\begin{thm} \label{t:eneq}
Let $a(t)$ be a solution of \eqref{Bsystem} for which $\|a(t)\|^3_{5/6}$ is
integrable on some interval $[T_1,T_2]$. Then $a(t)$ satisfies the energy equality
\begin{equation}
|a(t)|^2 = |a(t_0)|^2 + 2\int_{t_0}^t (f, a(\tau)) \, d\tau,
\end{equation}
for all $T_1\leq t_0 \leq t \leq T$.
\end{thm}
\begin{proof}
%Assume that there exists a solution $a(t)$ for which
%$\|a(t)\|^3_{1/3} \in L^1_{\mathrm{loc}}([0,\infty))$. Take any $T>0$. 

Let $a(t)$ be a solution satisfying the hypothesis of the theorem.
First, we recall the property of $l^p$ spaces which states that if
$p \geq q$ then   
$$ \|h\|_{l^p} \leq \|h\|_{l^q},$$
for all $h\in l^q$. 
We shall apply this property with $h_j = \lambda^{\frac{2j}{3}} a_j^2$
and $p=3/2$ to obtain 
\begin{equation} \label{ineq} 
\int_{t_0}^{t}\sum_{j=0}^\infty \lambda^{j} a_j^3 \, d\tau \leq
\int_{t_0}^{t}\left(\sum_{j=0}^\infty \lambda^{\frac{2j}{3}} a_j^2\right)^{3/2} \, d\tau < \infty,
\end{equation}
where $T_1\leq t_0 \leq t \leq T$.
However the expression \eqref{ineq} with $\lambda = 2^{5/2}$ combined with the assumption of the theorem 
implies that 
\begin{equation}
\int_{t_0}^t \lambda^{N} a_N^3 \, d\tau \to 0 \qquad \text{as} \qquad N \to \infty.
\end{equation}
Since
\begin{equation}
\int_{t_0}^t \lambda^{N}a^2_N a_{N+1} \, d \tau \leq
\int_{t_0}^t \lambda^{N}a^3_N \, d \tau + \int_{t_0}^t \lambda^{N+1}a^3_{N+1} \, d \tau,
\end{equation}
we can take the limit of \eqref{e:eipar} as $N\to \infty$ to obtain
\begin{equation}
|a(t)|^2 = |a(t_0)|^2 + 2\int_{t_0}^t (f, a(\tau)) \, d\tau.
\end{equation}
\end{proof}

As a consequence, we can now show that every solution (with any initial
data in $l^2$) blows up in finite time in $H^{5/6}$ norm.

\begin{cor}
Let $a(t)$ be a solution of \eqref{Bsystem} with $a_j(0)\geq 0$.
Then $\|a(t)\|^3_{5/6}$ is not
locally integrable on $[0, \infty)$. In particular, $\|a(t)\|_{5/6}$ blows up in finite time.
\end{cor}
\begin{proof}
Assume that $\|a(t)\|^3_{5/6}$ is locally integrable on $[0, \infty)$. Note that
Theorem~\ref{t:Expon}
implies that $a(t)$ converges to the fixed point in $l^2$. Therefore,
\begin{equation}
\lim_{t\to\infty}a_j(t) = \lambda^{-\frac{j}{3}}.
\end{equation}
In particular, we have that
\begin{equation}
\lim_{t\to\infty}a_0(t) =1.
\end{equation}
Hence, there exists $T>0$, such that $a_0(t) \geq 1/2$ for all $t\geq T$. Therefore,
\begin{equation}
(f,a(t)) \geq \frac{1}{2} \lambda^{-\frac{1}{3}}, \qquad \forall t\geq T.
\end{equation}
Hence, thanks to Theorem~\ref{t:eneq},
\begin{equation}
|a(t)|^2 \geq \frac{1}{2} \lambda^{-\frac{1}{3}}(t-T),   \qquad \forall t\geq T,
\end{equation}
which contradicts the fact that $a(t)$ converges to the fixed point in $l^2$
as $t \to \infty$
(Theorem~\ref{t:Expon}).
\end{proof}

\section{Global attractor} \label{secattr}
Since the uniqueness of a solution for given initial data is an open problem,
it is not known whether a semigroup of solution operators can be define for the
dyadic model. Therefore, we use a more general framework of an evolutionary system
$\Dc$ from \cite{CF,C2}.

%Let
%\begin{equation*}
%\ds(a,b)=\|a-b\|_{l^2}, \qquad
%\dw(a,b)= \sum_{j=0}^\infty \frac{1}{\lambda^{(n^2)}}
%\frac{\|a_j-b_j\|_{l^2}}{1 + \|a_j-b_j\|_{l^2}},
%\qquad a,b \in l^2.
%\end{equation*}
%Here, $\ds$ is a strong distance, and $\dw$ is a weak distance
%that induces a weak topology on any bounded subset of $l^2$.

Let $X$ be a closed ball in $H$ centered at the fixed point. Note that $X$ is
weakly compact.
In addition, 
\begin{equation} \label{e:inX}
a(t) \in X, \qquad \forall t\geq 0,
\end{equation}
for every solution $a(t)$ to \eqref{Bsystem} with the initial data
$a(0) \in X$.
%Moreover, we have the following.
%\begin{thm}
%$X$ is an absorbing ball, i.e., for any bounded subset $K \in l^2$, there
%exists a time $t_0\geq 0$, such that
%\begin{equation}
%a(t) \in X, \qquad \forall t\geq t_0,
%\end{equation} 
%for every solution $a(t)$ with $a(0) \in K$.
%\end{thm}
%\begin{proof}
%Assume the contrary. Then there exists a bounded set 
%\end{proof}

Let
\begin{equation}
C([0,\infty);\Xw):=\{a(\cdot): [0,\infty) \to X, a_n(t) \mbox{ continuous for all }
n\}
\end{equation}
endowed with the distance
\begin{equation}
d_{C([0, \infty);\Xw)}(a,b) = \sum_{T\in \mathbb{N}} \frac{1}{\lambda^T} \frac{\sup\{\dw(a(t),b(t)):0\leq t\leq T\}}
{1+\sup\{\dw(a(t),b(t)):0\leq t\leq T\}}.
\end{equation}
Let
\begin{equation}
\mathcal{T} := \{ I: \ I=[T,\infty) \subset \mathbb{R}, \mbox{ or } 
I=(-\infty, \infty) \},
\end{equation}
and for each $I \subset \mathcal{T}$ let $\mathcal{F}(I)$ denote
the set of all $X$-valued functions on $I$.
%Now we define an evolutionary system $\Dc$ as follows.
%\begin{definition} \label{Dc}
A map $\Dc$ that associates to each $I\in \mathcal{T}$ a subset
$\Dc(I) \subset \mathcal{F}$ will be called an evolutionary system if
the following conditions are satisfied:
\begin{enumerate}
\item $\Dc([0,\infty)) \ne \emptyset$.
\item
$\Dc(I+s)=\{a(\cdot): \ a(\cdot -s) \in \Dc(I) \}$ for
all $s \in \mathbb{R}$.
\item $\{a(\cdot)|_{I_2} : a(\cdot) \in \Dc(I_1)\}
\subset \Dc(I_2)$ for all
pairs of $I_1,I_2 \in \Omega$, such that $I_2 \subset I_1$.
%\item
%$H=\{ u(0): u(\cdot) \in \Dc([0, \infty)) \}$
\item
$\Dc((-\infty , \infty)) = \{a(\cdot) : \ a(\cdot)|_{[T,\infty)}
\in \Dc([T, \infty)) \ \forall T \in \mathbb{R} \}.$
\end{enumerate}
Let
\begin{equation*}
\begin{split}
R(t)A &:= \{a(t): a(0)\in A, a \in \Dc([0,\infty))\},\\
\Rc(t) A &:= \{a(t): a(0) \in A, a \in \Dc((-\infty,\infty))\}, \qquad A \subset X,
\ t \in \mathbb{R}.
\end{split}
\end{equation*}

A set $A \subset X$ uniformly attracts a set $B \subset X$
if for any $\epsilon>0$ there exists $t_0$, such that
\begin{equation}
R(t)B \subset B(A, \epsilon), \qquad \forall t \geq t_0.
\end{equation}

For $A \subset X$ and $r>0$ denote
$B_{\bullet}(A,r) = \{a: \ d_{\bullet}(A,a) < r\},$ where
$\bullet = \mathrm{s,w}$.
Now we define an attracting set and a global attractor as follows.

\begin{defn}
A set $A \subset X$ is $\mathrm{d}_{\bullet}$-attracting set
($\bullet = \mathrm{s,w}$) if it uniformly
attracts $X$ in $\mathrm{d}_{\bullet}$-metric, i.e.,
for any $\epsilon>0$ there exists
$t_0$, such that
\begin{equation}
R(t)X \subset B_{\bullet}(A, \epsilon), \qquad \forall t \geq t_0.
\end{equation}
A set $A \subset X$ is invariant if $\Rc(t) A = A$ for all $t\geq 0$.
A set $\mathcal{A}_{\bullet}\subset X$ is a
$\mathrm{d}_{\bullet}$-global attractor if
$\mathcal{A}_{\bullet}$ is a minimal $\mathrm{d}_{\bullet}$-closed
$\mathrm{d}_{\bullet}$-attracting  set.
%\item $\mathcal{A}_{\bullet}$ is $\mathrm{d}_{\bullet}$-closed.
%\end{enumerate}
\end{defn}
The following result was proved in \cite{CF}:
\begin{thm} \label{c:weakA}
The evolutionary system $\Dc$ always possesses a weak global attractor
$\Aw$.
In addition, if $\Dc([0,\infty))$ is compact in
$C([0, \infty);\Xw)$, then
\begin{enumerate}
\item[(a)]
$ \Aw =\{a^0: \ a^0=u(0) \mbox{ for some }
u \in \Dc((-\infty, \infty))\}.$
\item[(b)] $\Aw$ is the maximal invariant set.
\end{enumerate}
\end{thm}
%We will also use the following theorem from \cite{C}.
%\begin{thm}
%If $\Dc([0,\infty))$ is compact in $C([0, \infty);\Xw)$, the evolutionary system possesses
%an energy inequality, and every complete trajectory is strongly continuous, then
%$\Aw$ is a strongly compact strong global attractor.
%\end{thm}

For the dyadic model we define $\Dc$ in the following way.
\begin{equation*}
\begin{split}
\Dc([T,\infty)) &:= \{a: a(\cdot)
\mbox{ is a solution to \eqref{Bsystem} on } [T,\infty), a(0) \in X,
a_j(0)\geq 0\}, \\
\Dc((\infty,\infty)) &:=  \{a: a(\cdot)
\mbox{ is a solution to \eqref{Bsystem} on } (-\infty,\infty), a(0) \in X,
a_j(0)\geq0\}.
\end{split}
\end{equation*}
Clearly, $\Dc$ satisfies properties (1)--(4). Then Theorem~\ref{c:weakA} immediately 
yields that the weak global attractor $\Aw$ exists. In order to infer
that $\Aw$ is the maximal invariant set, we need the following result.

\begin{lem} \label{l:compact}
$\Dc([0,\infty))$ is compact in $C([0, \infty);\Xw)$.
\end{lem}
\begin{proof}
Take any sequence $a^k \in \Dc([0,\infty))$. 
Thanks to \eqref{e:inX}, there exists $R>0$, such that
\begin{equation}
a^{k}_j(t) \leq R, \qquad \forall n, t\geq 0.
\end{equation}
Therefore,
\begin{equation}
|a^k_j(t) - a^k_j(s)|
\leq (\lambda^{j-1} R^2 + \lambda^jR^2 + f_j)
|t-s|,
\end{equation}
for all $j$, $t\geq 0$, $s\geq0$. Thus,
\begin{equation}
\dw(a^{k}(t), a^{k}(s)) = \sum_{j=0}^\infty \frac{1}{\lambda^{(j^2)}}
\frac{|a^k_j(t)-a^k_j(s)|}{1 + |a^k_j(t)-a^k_j(s)|} \leq c|t-s|,
\end{equation}
for some constant $c$ independent of $k$. Hence, $\{a^{k}\}$ is an
equicontinuous sequence
of functions in $C([0,\infty);\Xw)$ with bounded initial data. Therefore,
Ascoli-Arzela theorem implies that
$\{a^k\}$ is relatively compact in $C([0,T]; \Xw)$, for all time $T>0$. Using
a diagonalization process, we obtain that $\{a^k\}$ is relatively compact in
$C([0,\infty); \Xw)$. Hence, there exists a weakly continuous $X$-valued
function $a(t)$ on $[0,\infty)$, such that
\begin{equation} \label{weakconv}
a^{k_n} \to a \qquad \mbox{as} \qquad k_n \to \infty  \qquad \mbox{in}
\qquad  C([0,\infty); \Xw),
\end{equation}
for some subsequence $k_n$.
In particular, since $X$ is weakly compact, $a(0) \in X$.

In addition, since $a^{k_n}(t)$ is a solution to \eqref{Bsystem},
we have
\begin{equation}
a^{k_n}_n(t) = a^{k_n}_j(0) + \int_0^t (
\lambda^{j-1} (a^{k_n}_{j-1})^2 - \lambda^{j} a^{k_n}_j a^{k_n}_{j+1} +f_j) \, d\tau,
\end{equation}
for all $j$. Taking the limit as $k_n \to \infty$, we obtain
\begin{equation}
a_j(t) = a_j(0) + \int_0^t (
\lambda^{j-1} a_{j-1}^2 - \lambda^{j} a_j a_{j+1} +f_j) \, d\tau,
\end{equation}
for all $j$.
Since $a_j(t)$ is continuous, $a_j \in C^1([0,\infty))$ and satisfies
\eqref{Bsystem}.
\end{proof}

Finally, we show that the weak global attractor is the fixed point.

\begin{thm}
A strong global attractor $\As$ of \eqref{Bsystem} exists, $\As=\Aw$,
and it is the fixed point: 
\begin{equation}
\As = \Aw = \{\hat{a}\},
\end{equation}
where $\hat{a}_j=\lambda^{-\frac{j}{3}}$. Moreover, $\As$ is an exponential attractor.
\end{thm}
\begin{proof}
The statement of the theorem immediately follows from Theorem~\ref{t:Expon}.
%Consider a complete trajectory $a \in \Dc((-\infty,\infty))$. Due to
%Theorem~\ref{decreasing-norm}, 
%\begin{equation*}
%\|b(t_2)\|^2_{l^2} -\|b(t_1)\|^2_{l^2} \leq
%- \alpha\int_{t_1}^{t_2} \|d(t)\|_{l^2}^2 \, dt,
%\end{equation*}
%Since $a(t) \in X$ for all $t \in \mathrm{R}$, we have that
%\begin{equation} \label{e:btoc} 
%\|b(t)\|_{l^2} \to c,\qquad \text{as} \qquad
%t \to -\infty,
%\end{equation}
%for some $c\geq 0$, and
%\begin{equation} \label{e:dto0}
%\|d(t)\|_{l^2} \to 0,\qquad \text{as} \qquad
%t \to -\infty.
%\end{equation}
%Therefore, thanks to \eqref{eq:bvsd}, we have that $\|b(t)\|_{l^2} \to 0$
%as $t \to - \infty$. Since $\|b(t)\|_{l^2}\|$ is decreasing, $b(t)=0$ for all
%$t\in\mathbb{R}$. Hence, $a(t) = \lambda^{-\frac{j}{3}}$ for all $t\in\mathbb{R}$.
\end{proof}

\section{Onsager's Conjecture and Kolmogorov's $5/3$ law} \label{secOns} 

As we discussed in \cite{CFP} much attention has been given to the 
statistical theories of turbulence developed by Kolmogorov \cite{K} and
Onsager \cite{Ons}. It is suggested that an appropriate mathematical
description of 3-dimensional turbulent flow is given by weak solutions
of the Euler equations which are not regular enough to conserve energy. Onsager 
conjectured that for the velocity H\"{o}lder exponent $h > 1/3$ the energy is
conserved and that this ceases to be true for $h \leq 1/3$. This phenomenon 
is now called turbulent or anomalous dissipation. Kolmogorov's theory 
predicts that in a fully developed turbulent flow the energy spectrum 
$E(|k|)$ in the inertial range is given by 
\be{KLrev} 
E(|k|) = c_0 \bar{\epsilon}^{2/3} |k|^{-5/3},
\end{equation} 
where $\bar{\epsilon}$ is the average of the energy dissipation rate. 
The model system that we study in this present paper exactly reproduces
the phenomena described above. The appropriate choice of $\lambda$ for the 
3-dimensional model is $2^{5/2}$. Interpreting the results of 
sections \ref{secfix} and \ref{secattr} with $\lambda = 2^{5/2}$, 
we proved that regular solutions, i.e.
solutions with bounded $H^{5/6}$ norm, satisfy the energy equality 
whereas solutions after the time of blow-up loose regularity and dissipate energy. 
The model \eqref{Bsystem} is derived under the assumption that $a_j^2(t)$ 
is the total energy in the shell $2^j \leq |k| < 2^{j+1}$. Hence
the energy spectrum for the fixed point is: 
\be{KLfix} 
E(|k|) = 2^{5/6} f_0 |k|^{-5/3}.
\end{equation}
Since the dissipation rate for the fixed point is equal to 
the energy input rate, we have 
$$ \bar{\epsilon} = a_0 f_0 = 2^{5/12} f_0^{3/2}.$$
Furthermore, since the support of any time average measure belongs to the global
attractor, the average dissipation rate is equal to the dissipation rate of the fixed 
point. This result is proved in \cite{FMRT} for the 3D Navier-Stokes equations,
but  it also holds for the inviscid dyadic model due to the anomalous dissipation.
Thus Kolmogorov's law is valid for the system.

\end{document}